\documentclass[a4paper,11pt]{amsart}

\usepackage{amsmath,amssymb,amsthm}
\usepackage{amsfonts}

\newtheorem{theorem}{Theorem}[section]

\newtheorem{proposition}[theorem]{Proposition}
\newtheorem{corollary}[theorem]{Corollary}
\newtheorem{e-definition}[theorem]{Definition\rm}
\def\remark{\nin {\bf Remark.\ \, }}

\def\ignore#1{}
\def\eq{\begin{equation}}
\def\en{\end{equation}}
\def\eqa{\begin{eqnarray}}
\def\ena{\end{eqnarray}}
\def\eqs{\begin{eqnarray*}}
\def\ens{\end{eqnarray*}}

\def\comp{{\mathbb C}}
\def\non{\nonumber}
\def\s{\sigma}
\def\l{\lambda}
\def\d{\delta}
\def\Ref#1{(\ref{#1})}
\def\Eq{\ =\ }
\def\Le{\ \le\ }

\def\h{\eta}

\def\th{\theta}
\def\g{\gamma}
\def\ch{\chi}
\def\k{\kappa}
\def\G{\Gamma}
\def\n{\nu}
\def\e{\varepsilon}
\def\f{\phi}
\def\a{\alpha}
\def\z{\zeta}

\def\nn{{\mathcal N}}
\def\tps{{\tilde\ps}}
\def\hps{{\hat\ps}}
\def\bps{{\bar\ps}}

\def\pr{{\mathbb P}}
\def\integ{{\mathbb Z}}
\def\ex{{\mathbb E}}
\def\ui{^{(1)}}
\def\ut{^{(2)}}

\def\dtv{d_{TV}}

\def\law{{\mathcal L}}

\def\r{\rho}
\def\D{\Delta}
\def\m{\mu}
\def\ps{\psi}

\def\Blb{\left\{}
\def\Brb{\right\}}
\def\Po{{\rm Po\,}}
\def\Be{{\rm Be\,}}

\def\nin{\noindent}
\def\half{{\textstyle\frac12}}

\def\Bl{\Bigl(}
\def\Br{\Bigr)}

\def\nat{{\mathbb N}}
\def\Def{\ :=\ }

\def\sso{\sum_{s\ge0}}
\def\sjo{\sum_{j\ge0}}

\def\slr{\sum_{l=1}^r}

\def\skol{\sum_{k=0}^l}

\def\un{^{(n)}}
\def\ta{\tilde a}
\def\tal{\tilde \a}

\def\var{{\rm Var\,}}

\def\p{\pi}

\def\ep{\hfill$\Box$\\[2ex]}

\def\nti{n\to\infty}

\def\b{\beta}

\def\pl{p_{\l}}
\def\zl{\z_{\l}}
\def\sjor{{\sum_{j=0}^r}}
\def\slor{{\sum_{l=0}^r}}

\def\impp{\int_{-\pi}^\pi}
\def\imii{\int_{-\infty}^\infty}
\def\pgf{probability generating function}
\def\andy{{\quad\mbox{and}\quad}}
\def\sji{\sum_{j\ge1}}

\def\CP{{\rm CP\,}}

\def\ur{^{(r)}}
\def\us{^{(s)}}
\def\ceka{\v Cekanavi\v cius}
\def\slt{\sum_{l\ge2}}

\def\dloc{d_{{\rm loc}}}
\def\dK{d_{{\rm K}}}
\def\BBB{\bar A}

\def\tA{{\widetilde A}}
\def\kkk{\g}
\def\mmm_#1{K_{#1}}

\def\slir{\sum_{l=1}^r}
\def\slor{\sum_{l=0}^r}
\def\loglog{\log\log}
\def\Blm{\left|}
\def\Brm{\right|}
\def\tC{{\widetilde C}}
\def\om{\omega}
\def\Om{\Omega}
\def\lnti{\lim_{n\to\infty}}
\def\dmn{d_{\m\n}}

\begin{document}

\title{Mod--discrete expansions}

\author{A.\ D.\ Barbour} \address{Institut f\"ur Mathematik,
  Universit\"at Z\"urich\\Winterthurertrasse 190, CH-8057 Z\"URICH,
  Switzerland} \email{a.d.barbour@math.uzh.ch}

\author{E.\ Kowalski} \address{ETH Zurich, D-MATH\\ R\"amistrasse 101,
  8092 Z\"URICH, Switzerland} \email{kowalski@math.ethz.ch}
\thanks{Work of EK supported in part by the National Science
  Foundation under agreement No. DMS-0635607 during a sabbatical stay
  at the Institute for Advanced Study}

\author{A. Nikeghbali} \address{Institut f\"ur Mathematik,
  Universit\"at Z\"urich\\Winterthurertrasse 190, CH-8057 Z\"URICH,
  Switzerland} \email{ashkan.nikeghbali@math.uzh.ch} \thanks{Work of
  ADB and AN supported in part by Schweizerischer Nationalfonds
  Projekte Nr.\ 20--117625/1 (ADB) and Schweizerischer Nationalfonds
  Projekte Nr.\ 200021\_119970/1 (AN)}


\date{}

\subjclass[2000]{62E17; 60F05, 60C05, 60E10, 11N60}

\keywords{mod--Poisson convergence; characteristic function;
  Poisson--Charlier expansion; Erd\H os--Kac theorem} 

\begin{abstract}
  In this paper, we consider approximating expansions for the
  distribution of integer valued random variables, in circumstances in
  which convergence in law cannot be expected. The setting is one in
  which the simplest approximation to the $n$'th random variable~$X_n$
  is by a particular member $R_n$ of a given family of distributions,
  whose variance increases with~$n$.  The basic assumption is that the
  ratio of the characteristic function of~$X_n$ and that of~$R_n$
  converges to a limit in a prescribed fashion.  Our results cover a
  number of classical examples in probability theory, combinatorics
  and number theory.  
\end{abstract}


\maketitle

\setcounter{equation}{0}
\section{Introduction}

In a remarkable paper, Hwang~(1999) considered sequences
of non-nega\-tive integer valued random variables~$X_n$, whose probability
generating functions~$f_{X_n}$ satisfy
\[
      e^{\l_n(1-z)}f_{X_n}(z) \ \to \ g(z), 
\]
for all $z\in\comp$ with $ |z| \le \h > 1$, where the function~$g$ is
analytic, and $\lnti\l_n = \infty$.  Under some extra conditions, he
exhibits tight bounds on the accuracy of the approximation of the
distribution of~$X_n$ by a Poisson distribution with carefully chosen
mean, close to~$\l_n$. Independently, motivated by specific examples
arising in Random Matrix Theory and number theory, Jacod, Kowalski and
Nikeghbali~(2008) explored the properties of a related ratio
convergence for real valued random variables, namely when the
characteristic functions~$\f_{X_n}$ satisfy
\[
     e^{-i\th\b_n + \th^2\g_n/2}\f_{X_n}(\th) \ \to\ \Phi(\th),
\]
locally uniformly in~$\th$ (in particular, bounds on the error in the
approximation of the distribution~$P_{X_n}$ by the normal distribution
$\nn(\b_n,\g_n)$ can be simply deduced).  Kowalski and
Nikeghbali~(2009) went on to explore some consequences (and structural
aspects in arithmetic cases) of the corresponding uniform limit
\eq\label{mod-Poisson-cvgce} \exp\{\l_n(1-e^{i\th})\}\f_{X_n}(\th) \to
\ps(\th), \quad 0 < |\th| \le \pi, \en for random variables~$X_n$
(usually integer valued), as in Hwang~(1999) with the Poisson
characteristic function in the ratio.  Note that the conditions on the
distributions of the~$X_n$ are now much weaker than those of
Hwang~(1999): for instance, his conditions require the~$X_n$ to take
only non-negative values, and to have exponential tails.  On the other
hand, the probabilistic results that Hwang derives are much more
sophisticated.  He establishes bounds on the error in his
approximations with respect to a number of different metrics, and
shows that they are sharp.  For instance, for the Kolmogorov and total
variation distances, his bounds are typically of order~$O(\l_n^{-1})$,
and he also gives the value of the leading asymptotic term in the
error.

In this paper, we work with integer valued random variables, and with 
characteristic function conditions that sharpen~\Ref{mod-Poisson-cvgce},
with the aim of developing approximations of higher order.  Our main
result, Proposition~\ref{Th0}, is very simple and explicit.  This
enables us to dispense with asymptotic settings, and to prove concrete
error bounds. As a direct consequence, we are able to deduce a Poisson--Charlier
approximation with error of order $O(\l_n^{-(r+1)/2})$, for any 
prescribed~$r$, assuming that Hwang's conditions hold. Our Poisson--Charlier 
expansions are
derived under more general conditions, in which the~$X_n$ may
have only a few finite moments. These are established in Section~\ref{PC}, and simpler,
translated Poisson approximations are considered in Section~\ref{Poisson}.

Hwang~(1999) notes that his methods are also applicable to families of
distributions other than the Poisson family, and gives examples using
the Bessel family.  Our approach allows one to derive expansions based
on any discrete family of distributions, as shown in Section~\ref{CP},
provided that their characteristic functions satisfy a simple
condition, and this without any extra effort.  Indeed, the main
problem is to identify the higher order terms in the expansions.
These turn out to be simply the higher order differences of the basic
distribution, leading, for example, to the Charlier polynomial factors
in the Poisson case.  We discuss some examples, to sums of independent
integer valued random variables, to Hwang's setting and to the 
Erd\H os--Kac theorem, in Section~\ref{Apps}.

\medskip
\remark We recall the motivation behind the terminology (mod-gaussian,
mod-poisson, and here mod-discrete): the simplest example leading to
limits like (say)~(\ref{mod-Poisson-cvgce}) is when $X_n=P_n+Y$ where
$P_n$ has Poisson distribution $\Po(\l_n)$ and is independent of $Y$,
where $\ps(\th)$ is the characteristic function of $Y$. Thus the
sequence converges to $Y$ ``modulo Poisson variables''.

\section{The basic estimate}
\setcounter{equation}{0}

We frame our approximations in terms of three distances between (signed)
measures $\m$ and~$\n$ on the integers: the point metric
\[
    \dloc(\m,\n) \Def \sup_{j\in\integ}|\m\{j\} - \n\{j\}|,
\]
the Kolmogorov distance
\[
    \dK(\m,\n) \Def \sup_{j\in\integ}|\m\{(-\infty,j]\} - \n\{(-\infty,j]\}|,
\]
and the total variation norm
\[
    \|\m-\n\| \Def \sum_{j\in\integ}|\m\{j\} - \n\{j\}|
=2\sup_{A\subset \integ}{|\m(A)-\n(A)|}.
\]
Other metrics could also be treated using our results.  Our conditions
are expressed in terms of characteristic functions, defined,
for a finite signed measure~$\s$ on~$\integ$, by 
$\f_\s(\th) := \sum_{j\in\integ} e^{ij\th}\s\{j\}$, for
$|\th|\le\p$.  The essence of our argument is the following simple result, 
linking
the closeness of the signed measures to the closeness of their
characteristic functions, when these have a common factor involving
a `large' parameter~$\r$.

\begin{proposition}\label{Th0}
Let $\m$  and $\n$ be finite signed measures on~$\integ$, with
characteristic functions $\f_\m$ and~$\f_\n$ respectively. 
Suppose that $\f_\m = \ps_\m \ch$ and $\f_\n = \ps_\n \ch$,
where, for some $ \kkk_1, \kkk_2,\r,t>0$,
\[
     |\ps_\m(\th) - \ps_\n(\th)| \Le  \kkk_1|\th|^t \quad\mbox{and}\quad
     |\ch(\th)| \le  \kkk_2 e^{-\r\th^2} \quad \mbox{for all}\quad |\th|\le\p.
\]
Then, writing $ \kkk = \kkk_1 \kkk_2$, there are explicit constants
$\a_{1t}$ and $\a_{2t}$ such that \eqs
1.&& \sup_{j\in \integ}|\m\{j\} - \n\{j\}| \Le \a_{1t}  \kkk (\r\vee1)^{-(t+1)/2};\\
2.&& \sup_{a \le b\in \integ}|\m\{[a,b]\} - \n\{[a,b]\}| \Le \a_{2t}
\kkk (\r\vee1)^{-t/2}.  \ens
\end{proposition}

\proof
For any $j\in\integ$, the Fourier inversion formula gives
\eq\label{Fourier-1}
   \m\{j\} - \n\{j\} 
     \Eq \frac1{2\pi}\impp e^{-ij\th}(\ps_\m(\th)-\ps_\n(\th))\ch(\th)\,d\th,
\en			
from which our assumptions imply directly that
\eqs
  |\m\{j\} - \n\{j\}| &\le& 
	    \frac1{2\pi}\impp  \kkk |\th|^t\exp\{-\r\th^2\}\,d\th.
\ens
For $\r\le1$, we thus have
\eqs
   |\m\{j\} - \n\{j\}| &\le& 
	    \frac1{2\pi}\impp  \kkk |\th|^t\,d\th 
   \Le \frac{\pi^{t} \kkk}{t+1} \ =:\ \b_{1t} \kkk. 
\ens
For $\r\ge1$, it is immediate that
\eqs
  |\m\{j\} - \n\{j\}| &\le& \frac{ \kkk}{2\pi}\Bl \frac{1}{\sqrt{2\r}}\Br^{t+1}
	    \imii |y|^t e^{-y^2/2}\,dy 	\Le \b'_{1t}  \kkk \r^{-(t+1)/2},
\ens
with $\b'_{1t} := 2^{-(t+1)/2}m_{t}/\sqrt{2\pi}$;  here, $m_t$ denotes the
$t$-th absolute moment of the standard normal distribution. 
Setting $\a_{1t} := \max\{\b_{1t},\b'_{1t}\}$, this proves part~1.
The second part is similar, adding~\Ref{Fourier-1} over $a\le j\le b$, and estimating
\[
    \frac{\bigl|e^{-ia\th} - e^{-i(b+1)\th}\bigr|}{|1 - e^{-i\th}|} \Le
         \frac{\pi}{|\th|},
		\qquad |\th| \le \pi.
\]
This gives part~2, with 
\[
   \a_{2t} := \max\{2^{-t/2}m_{t-1}\sqrt{\pi/2}, \pi^{t}/t\}.
\]

\bigskip
\nin In particular, the second part bounds the distance between the two
measures in the Kolmogorov distance.  We shall principally be concerned
with taking~$\m$ to be the distribution of a random variable~$X$; we
allow~$\n$ to be a signed measure largely for reasons of technical 
convenience.

For some applications, a slight weakening of the conditions in
Proposition~\ref{Th0} is useful.  The following result is
proved in exactly the same way as before.

\medskip
\begin{proposition}\label{Th0'}
Let $\m$  and $\n$ be finite signed measures on~$\integ$, with
characteristic functions $\f_\m$ and~$\f_\n$ respectively. 
Suppose that $\f_\m = \ps_\m \ch$ and $\f_\n = \ps_\n \ch$,
where, for some $\th_0, \kkk,\e,\h,\r'>0$ and for 
positive pairs $\kkk_m,t_m$, $1\le m\le M$, we have
\begin{align*}
     |\ps_\m(\th) - \ps_\n(\th)| &\le  \sum_{m=1}^M \kkk_m|\th|^{t_m} 
       + \e\quad\mbox{and}\quad |\ch(\th)| \le  \kkk e^{-\r\th^2},
    \quad  0 \le |\th|\le \th_0;\\
    |\f_\m(\th) - \f_\n(\th)| &\le \h, \qquad \th_0 < |\th| \le \pi.
  \end{align*}
Then, with notation as for Proposition~\ref{Th0}, we have
\eqs
  1.&& \sup_{j\in \integ}|\m\{j\} - \n\{j\}| 
     \Le \sum_{m=1}^M \kkk_m \kkk\a_{1t_m}  (\r\vee1)^{-(t_m+1)/2} 
      + \tal_1  \kkk\e + \tal_2 \h;\\[2ex]
  2.&& \sup_{a_0 \le a \le b \le b_0}|\m\{[a,b]\} - \n\{[a,b]\}| \\
    &&\hskip0.2in \Le \sum_{m=1}^M \kkk_m \kkk\a_{2t_m}  (\r\vee1)^{-t_m/2}
      + (b_0-a_0+1)(\tal_1  \kkk\e + \tal_2 \h),
\ens
where
\[
   \tal_{1} \Def \Bl \frac{\th_0}{\pi} \wedge \frac1{2\sqrt{\pi\r}} \Br;
   \qquad \tal_2 \Def \Bl 1 - \frac{\th_0}{\pi} \Br.
\]
\end{proposition}

\bigskip
The presence of the factor $(b_0-a_0+1)$ in the second bound means that
a direct bound on the Kolmogorov distance between the signed measures
$\m$ and~$\n$ is not immediately visible.  
The following corollary is however easily deduced.

\medskip
\begin{corollary}\label{Cor0}
Under the conditions of Proposition~\ref{Th0'}, 
\eqs
   \dK(\m,\n) &\le& \inf_{a\le b}\bigl( \e^{{\rm (K)}}_{ab} + |\m|\{(-\infty,a)\cup(b,\infty)\}
     + |\n|\{(-\infty,a)\cup(b,\infty)\}\bigr); \\
   \|\m-\n\| &\le& \inf_{a\le b}\bigl( \e\ui_{ab} + |\m|\{(-\infty,a)\cup(b,\infty)\}
     + |\n|\{(-\infty,a)\cup(b,\infty)\}\bigr),
\ens
where 
\eqs
   \e^{{\rm (K)}}_{ab} &:=& \sum_{m=1}^M\kkk_m \kkk\a_{2t_m}  (\r\vee1)^{-t_m/2}
        + (b-a+1)(\tal_1  \kkk\e + \tal_2 \h);\\
   \e\ui_{ab} &:=& (b-a+1)\Blb\sum_{m=1}^M\kkk_m \kkk\a_{1t_m}(\r\vee1)^{-(t_m+1)/2}
        + (\tal_1  \kkk\e + \tal_2 \h)\Brb.
\ens
If also~$\m$ is a probability measure, then
\[
   \dK(\m,\n) \Le \inf_{a\le b}\bigl(1 - \n\{[a,b]\} + 2\e^{{\rm (K)}}_{ab} + 
     |\n|\{(-\infty,a)\cup(b,\infty)\}\bigr).
\]
\end{corollary}

\proof
The inequality for the total variation norm is immediate.  For the
Kolmogorov distance,
by considering the possible positions of $x$ in relation to $a < b$,
we have
\eqs
   \lefteqn{|\m\{(-\infty,x]\} - \n\{(-\infty,x]\}|}\\
  &&\Le \sup_{y < a}|\m\{(-\infty,y]\} - \n\{(-\infty,y]\}|
  + \sup_{a \le y \le b}|\m\{[a,y]\} - \n\{[a,y]\}|\\
   &&\qquad\quad\mbox{}
   + \sup_{y > b}|\m\{(b,y]\} - \n\{(b,y]\}| \\
  &&\Le |\m|\{(-\infty,a)\cup(b,\infty)\} + |\n|\{(-\infty,a)\cup(b,\infty)\}
    + \e^{{\rm (K)}}_{ab}.
\ens
If~$\m$ is a probability measure, we have 
\eqs
|\m|\{(-\infty,a)\cup(b,\infty)\} &=& 1 - \m\{[a,b]\} 
      \le 1 - \n\{[a,b]\} + \e^{{\rm (K)}}_{ab} .  \hskip1.5cm\Box
\ens

\bigskip
Under slightly stronger conditions than those of Proposition~\ref{Th0}, 
a much neater total variation bound can be deduced.

\medskip
\begin{proposition}\label{Th0-tv}
Let $\m$  and $\n$ be finite signed measures on~$\integ$, with
characteristic functions $\f_\m = \ps_\m \ch$ and~$\f_\n = \ps_\n \ch$ respectively,
where $\ch(\th) := \kkk_2 e^{-u(\th)}$ for some $\kkk_2 > 0$, and $u(0)=0$. 
Suppose now that~$u$ and the difference $\dmn := \ps_\m - \ps_\n $ are both twice
differentiable, that $u'(0) =  \dmn'(0) = 0$ and that, for some 
$ \kkk_1, \kkk_2, \kkk_3>0$, $\r\ge1$ and $t\ge2$,
\[
     |\dmn''(\th)| \Le  \kkk_1|\th|^{t-2},\quad |u''(\th)| \Le \kkk_3\r 
        \quad\mbox{and}\quad
     u(\th) \ge  \r\th^2 \quad \mbox{for all}\quad |\th|\le\p.
\]
Then, writing $ \kkk =  \kkk_1 \kkk_2$, there is a constant $\a_{3}:=\a_{3}(t,\kkk_3)$ such that
\eqs
  \|\m-\n\| \Le \a_{3}  \kkk (\r\vee1)^{-t/2}.
\ens
\end{proposition}

\medskip
\proof
First, the assumptions on~$\dmn$ and~$u$ give
\eq\label{ADB-new}
\begin{split}
   |\dmn'(\th)| &\le \frac{\kkk_1}{t-1}\,|\th|^{t-1};
    \quad |\dmn(\th)| \Le \frac{\kkk_1}{t(t-1)}\,|\th|^{t};\\
     |u'(\th)| &\Le  \kkk_3\r |\th|.
  \end{split}
  \en
In particular, for $|j| \le \lceil \sqrt\r \rceil$, we can apply
part~1 of Proposition~\ref{Th0}, which gives
\eq\label{near-bnd}
    |\m\{j\} - \n\{j\}| \Le \frac{\a_{1t}  \kkk}{t(t-1)} (\r\vee1)^{-(t+1)/2}.
\en
For the remaining~$j$, integrating the Fourier inversion 
formula~\Ref{Fourier-1} twice by parts gives
\begin{multline}
  \m\{j\} - \n\{j\} 
     = -\frac1{2\pi j^2}\impp e^{-ij\th}\Bigl(\dmn''(\th) - 2\dmn'(\th)u'(\th)
+\\
 \dmn(\th)\{(u'(\th))^2 - u''(\th)\}\Bigr)\ch(\th)\,d\th.
\label{Fourier-1-tv}
\end{multline}
Substituting the bounds from~\Ref{ADB-new} into~\Ref{Fourier-1-tv} gives	
\eqs
   \lefteqn{|\m\{j\} - \n\{j\}|} \\
    &\le&  \frac1{2\pi j^2}\impp \kkk 
         \Blb |\th|^{t-2} + \frac{2\kkk_3\r}{t-1}\,|\th|^t
             + \frac{\kkk_3\r}{t(t-1)}\,|\th|^t(1 + \kkk_3\r\th^2)\Brb
           e^{-\r\th^2}\,d\th  \\
    &\le& \frac\r{j^2}\,\kkk\,\b_3(t,\kkk_3)\r^{-(t+1)/2},
\ens
after some calculation, where, with $m_t$ as in Proposition~\ref{Th0},
\[
      \b_3(t,\kkk_3) \Def \frac{m_{t-2}}{4t\,2^{t/2}\sqrt\pi}\{4t+(2t+1)\kkk_3
         + (t+1)\kkk_3^2\}.
\]
Hence
\[
     \sum_{|j| > \lceil\sqrt\r\rceil} |\m\{j\} - \n\{j\}| \Le
           2\kkk\,\b_3(t,\kkk_3)\r^{-t/2},       
\]
and the proposition follows directly, with 
$\a_3(t,\kkk_3) := 2\b_3(t,\kkk_3) + \frac{5\a_{1t}}{t(t-1)}$.
\ep

\setcounter{equation}{0}
\section{Poisson--Charlier expansions}\label{PC}
Suppose first that~$X$ is an integer valued random variable having
characteristic function~$\f_X$ of the form $\f_X(\th) = \ps(\th)\pl(\th)$,
$|\th| \le \pi$, where~$\pl$ denotes the characteristic function of the Poisson
distribution~$\Po(\l)$ with mean~$\l$.  Underlying our considerations
is an unspecified asymptotic setting in which~$\l$
is large and~$\ps$ is thought of as (almost) fixed, but we do not need
to make direct use of this.  We now assume in addition that, for some $r\in\nat_0$
and for some $K_{r\d} > 0$, $0 < \d \le 1$,
\eq\label{local-bnd}
   |\ps(\th) - \ps_r(\th)| \Le K_{r\d} |\th|^{r+\d},\qquad |\th| \le \pi,
\en	 
where
\eq\label{psr-def}
  \ps_r(\th) \Def \slor a_l(i\th)^l
\en
is a polynomial of degree~$r$ with real coefficients~$a_l$, thus implying that 
\hbox{$a_0=1$.}
If~$\ps$ is itself the characteristic function of a probability measure,	
this assumption roughly corresponds to assuming that the measure has (at least)~$r$ 
finite moments.  Alternatively, we could assume that
\eq\label{local-bnd-tilde}
   |\ps(\th) - \tps_r(\th)| \Le K_{r\d} |\th|^{r+\d},\qquad |\th| \le \pi,
\en	 
where
\eq\label{psr-def-tilde}
  \tps_r(\th) \Def \slor \ta_l\bigl(e^{i\th}-1\bigr)^l,
\en
again with real coefficients~$\ta_l$ and $\ta_0=1$.  If~$r=0$, and thus 
$\ps_0(\th) = 1$ for all~$\th$, we could now immediately use~\Ref{local-bnd} 
in conjunction with Proposition~\ref{Th0} to
approximate the distribution of~$X$ by the Poisson distribution~$\Po(\l)$,
with an error in Kolmogorov distance of order~$\l^{-\d/2}$; note that
\eq\label{rho-def}
    |\pl(\th)| \Eq \exp\{-\l(1-\cos\th)\} \Le e^{-\r\th^2},\quad |\th| \le\p,
\en
with $\r := 2\p^{-2}\l$.

We now want to go further, and use~\Ref{local-bnd} with higher values
of~$r$ to justify more sophisticated approximations with a higher order
of accuracy.  In order to do so, we need to find `nice' signed 
measures~$\n_r$, whose characteristic functions are 
at least as close to $\ps_r(\th)\pl(\th)$
and $\tps_r(\th)\pl(\th)$ as~$\ps(\th)\pl(\th)$ is.
Now $\tps_r(\th)\pl(\th)$ is itself the characteristic function of
a signed measure, which we can then take as our choice of~$\n_r$. To identify~$\n_r$,
observe that,
if~$\f_\m$ is the characteristic function of a
signed measure~$\m$, then $\bigl(e^{i\th}-1\bigr)^l \f_\m(\th)$ is the
characteristic function of the $l$'th difference $\D^l\m$ of~$\m$: 
\[
     \D^l\m\{j\} \Def \sum_{k=0}^l {l\choose k}(-1)^k \m\{j-l+k\}.
\]
For $\m$ the Poisson distribution, this yields the Poisson--Charlier 
signed measures: 
\eq\label{nu-char}
    \tps_r(\th)\pl(\th) \Eq \slor \ta_l\bigl(e^{i\th}-1\bigr)^l \pl(\th)
\en
is the characteristic function
of the signed measure $\n = \n_r(\l;\ta_1,\ldots,\ta_r)$ 
on $\nat_0$ defined by
\eq\label{nu-def}
  \n\{j\} \Def \Po(\l)\{j\}\Bigl\{1 + \slr (-1)^l\ta_lC_l(j;\l) \Bigr\},
\en
where 
\eq\label{Charlier-def}
    C_l(j;\l) \Def \skol (-1)^k{l\choose k}{j\choose k}k!\,\l^{-k} 
\en		
denotes the $l$-th Charlier polynomial (Chihara~1978, (1.9), p.~171).

Note that, if ${j\choose k}$ is replaced by $j^k/k!$ in~\Ref{Charlier-def},
one obtains the binomial expansion of $(1-j/\l)^l$.
As this suggests,
the values of $C_l(j;\l)$ are in fact small for~$j$ near~$\l$ if~$\l$ is large: 
\eq\label{BC6.1}
   |C_l(j;\l)| \Le 2^{l-1}\{|1-j/\l|^l + (l/\sqrt\l)^l\} 
\en
(Barbour \& \v Cekanavi\v cius 2002, Lemma~6.1). \Ref{BC6.1} thus implies that 
the $l$-th term in the sum in~\Ref{nu-def} has total variation norm
at most $|\ta_l|c_l\l^{-l/2}$, for a universal constant~$c_l$.  It also
implies that, in any interval of the form $|j-\l| \le c\sqrt\l$, which is 
where the probability mass of~$\Po(\l)$ is mostly to be found, the correction 
to the Poisson measure~$\Po(\l)$ is of uniform relative order~$O(\l^{-l/2})$. 
Indeed, the Chernoff inequalities for $Z \sim \Po(\l)$ can be expressed in 
the form
\eq
\max\{\pr[Z > \l(1+\d)], \pr[Z < \l(1-\d)]\}
	\Le \exp\{-\l\d^2/2(1+\d/3)\}, \label{Chernoff}
\en 
for $0 < \d \le 1$ (Chung \&~Lu 2006, Theorem~3.2).  Since also, from~\Ref{Charlier-def},
\[
     |C(j;\l)| \Le (1+j/\l)^l \Le 2^l \quad\mbox{if}\quad 0\le j\le \l,
\]
and since
\[
   {j\choose k}k!\,\l^{-k}\,\frac{e^{-\l}\l^j}{j!} \Eq \frac{e^{-\l}\l^{j-k}}{(j-k)!}
     \Le \frac{e^{-\l}\l^{j-l}}{(j-l)!}
\]
if $0\le k\le l$ and $j \ge l + \l$, it follows that, for any $l\ge0$,
we have
\begin{align*}
   \sum_{j=0}^m |C_l(j;\l)|\Po(\l)\{j\} &\le 2^l \pr[Z \le m]
	     \Le 2^l \exp\{-(\l-m)^2/3\l\}\\
\intertext{for $m\le \l,$ and}
\sum_{j\ge m} |C_l(j;\l)|	\Po(\l)\{j\} &\le 2^l \pr[Z \ge m-l]	  
        \Le 2^l \exp\{-(m-l-\l)^2/3\l\},
\end{align*}
for $\l + l \le m \le 2\l + l$.
\par
Writing $|\n|$ to denote the absolute measure associated with~$\n$, it
thus follows that
\eqa
  |\n|\{[0,m]\} &\le& \BBB_r\, e^{-(\l-m)^2/3\l},\quad\qquad 0\le m\le \l; \non\\
	|\n|\{[m,\infty)\} &\le& \BBB_r\, e^{-(m-r-\l)^2/3\l},\qquad \l+r\le m\le 2\l,
	   \label{tails-of-nu}
\ena
where $\BBB_r := 1 + \slr 2^l|\ta_l|$, demonstrating concentration of measure
for~$\n$ on a scale of~$\sqrt\l$ around~$\l$.  Moreover, it can be deduced	
from~\Ref{BC6.1} that there exists a positive constant 
$d = d(\ta_1,\ldots,\ta_r)$  such that
$\n\{j\} \ge 0$ for $|j-\l| \le d\l$, and it follows from~\Ref{tails-of-nu} that
$|\n|\{j\colon\,|j-\l| > d\l\} = O(e^{-\a\l})$ for some $\a > 0$.  Since also
$\n\{\nat_0\} = 1$, it thus follows that, even if~$\n$ is formally a signed measure, 
it differs from a probability only on a set of measure exponentially small
with~$\l$.	

Thus, if~\Ref{local-bnd-tilde} holds, it follows that 
$X$ has characteristic function $\ps(\th)\pl(\th)$ and 
$\n := \n_r(\l;\ta_1,\ldots,\ta_r)$ has characteristic function
$\f_\n = \tps_r(\th)\pl(\th)$, and that the conditions of 
Proposition~\ref{Th0} are satisfied with $\m = P_X$, $t=r+\d$, $\k = K_{r\d}$
and $\r = 2\p^{-2}\l$, this last from~\Ref{rho-def}.
If, instead, we are given the inequality~\Ref{local-bnd}, we can
write $e^{i\th}-1 = i\th \sso (i\th)^s/(s+1)!$, and equate the 
coefficients of $(i\th)^j$ in~\Ref{psr-def} with those for $1\le j\le r$ 
in~\Ref{nu-char}, giving  $\ta_1,\ldots,\ta_r$ implicitly in terms of $a_1,\ldots,a_r$:
\eq\label{b-def}
   a_j \Eq \sum_{l=1}^j \ta_l\sum_{(s_1,\ldots,s_l)\in S_{j-l}} \prod_{t=1}^l
        \frac1{(s_t+1)!},
\en
where $S_m := \bigl\{(s_1,\ldots,s_l)\colon\,\sum_{t=1}^l s_t = m\bigr\}$.
With this choice of $\ta_1,\ldots,\ta_r$, it follows that	
$\n = \n_r(\l;\ta_1,\ldots,\ta_r)$ has characteristic function~$\f_\n$ satisfying
\eq\label{Gr}
   |\ps_r(\th) - \f_\n(\th)| \Le \G_r |\th|^{r+1},\qquad |\th| \le \pi,
\en
for $\G_r := \G_r(a_1,\ldots,a_r)$.  Hence, in this case, we obtain 
\eq\label{chf-approx}
   |\ps(\th) - \f_\n(\th)| \Le (K_{r\d}+G_{r\d}) |\th|^{r+\d},\qquad |\th| \le \pi,
\en
with $G_{r\d} := \G_r\p^{1-\d}$, and the conditions of 
Proposition~\ref{Th0} are satisfied with $\m = P_X$, $t=r+\d$, $ \kkk = K_{r\d}+G_{r\d}$
and $\r = 2\p^{-2}\l$.
Thus, if either \Ref{local-bnd} or~\Ref{local-bnd-tilde} is satisfied, 
a signed measure from the family $\n_r(\l;b_1,\ldots,b_r)$ can be found,
which approximates the probability measure~$P_X$ in the sense implied
by Proposition~\ref{Th0}.  These measures are themselves rather explicit
perturbations of the Poisson distribution~$\Po(\l)$.

We summarize these considerations in the following theorem, which is
deduced directly from Proposition~\ref{Th0}.  Note that we 
shall later be primarily concerned
with applications in which $K_{r\d}$ and~$G_{r\d}$ are not small, and in which 
therefore~$\l$ must be big, if our bounds are to be useful. However, for
the sake of completeness, we phrase our bounds in a form which also allows
for accuracy of approximation if $K_{r\d}+G_{r\d}$ is small.		

\medskip
\begin{theorem}\label{Th1}
Let $X$ be a random variable on~$\integ$ with distribution~$P_X$.  Suppose that 
its characteristic function $\f_X$ is of the form~$\ps\pl$, where 
$\pl(\th)$ is the characteristic
function of the Poisson distribution~$\Po(\l)$ with mean~$\l$. Suppose also 
that \Ref{local-bnd} is satisfied, for some $r\in\nat_0$ and $\d\ge0$.
Let $\n_r := \n_r = \n_r(\l;\ta_1,\ldots,\ta_r)$ be as in~\Ref{nu-def},
with $\ta_1$,\ldots, $\ta_r$
given implicitly by~\Ref{b-def}. Then, writing $t=r+\d$, we have 
\eqs
  1.&& \dloc(P_X,\n_r) \Def \sup_{j\in \integ}|P_X\{j\} - \n_r\{j\}|
  \\
&& \mbox{}\hskip0.8in \Le \a'_{1t}
      (K_{r\d}+G_{r\d})(\l\vee1)^{-(t+1)/2};\\
  2.&& \dK(P_X,\n_r) \Def \sup_{l\in \integ}|P_X\{(-\infty,l]\} - \n_r\{[0,l]\}| \\
    &&\mbox{}\hskip0.8in \Le \a'_{2t} (K_{r\d}+G_{r\d}) (\l\vee1)^{-t/2},
\ens
with 
\begin{gather*}
\a'_{1t} := \a_{1t}(\p^2/2)^{(t+1)/2},\quad \a'_{2t} = \a_{2t} 
(\p^2/2)^{t/2},
\\
G_{r\d} := \G(a_1,\ldots,a_r)\pi^{1-\d}.
\end{gather*}
If~\Ref{local-bnd} is replaced by~\Ref{local-bnd-tilde}, the corresponding bounds
hold with $G_{r\d}=0$. 
\end{theorem}

\medskip
Theorem~\ref{Th1} enables one to deduce simple bounds for other measures of the
distance between $P_X$ and~$\n$.  For instance, for the total variation norm,
with judicious choice of $m_1$ and~$m_2$, we can use part~1 to bound
\eq\label{TV-1}
   \sum_{j=m_1+1}^{m_2-1}|P_X\{j\} - \n\{j\}| 
       \Le (m_2-m_1-1)\sup_{j\in\nat_0} |P_X\{j\} - \n\{j\}|,
\en
and then~\Ref{tails-of-nu} and part~2 to take care of the remaining tail probabilities:
\begin{align}
  \sum_{j\le m_1}|P_X\{j\} - &\n\{j\}| \le P_X\{(-\infty,m_1]\} + |\n|\{[0,m_1]\}
       \label{TV-2a}\\
&    \le  \sup_{l\in \nat_0}|P_X\{(-\infty,l]\} - \n\{[0,l]\}| + 2|\n|\{[0,m_1]\},
    \nonumber
  \end{align}
and
\begin{align}
   \sum_{j\ge m_2}|P_X\{j\} - &\n\{j\}| \le P_X\{[m_2,\infty)\} +
       |\n|\{[m_2,\infty)\}
        \label{TV-2}\\
    &\le  \sup_{l\in \nat_0}|P_X\{(-\infty,l]\} - \n\{[0,l]\}| + 2|\n|\{[m_2,\infty)\}.
    \nonumber
  \end{align}
This gives the following theorem.

\medskip
\begin{theorem}\label{Th2}
Suppose that the conditions of Theorem~\ref{Th1} are satisfied,
with~\Ref{chf-approx} holding. If $K_{r\d}+G_{r\d} < 1$, there is a 
constant~$\a_{4t}$ such that
\begin{multline}\label{Th2-1}
   \|P_X-\n\|
  \Le \a_{4t}(K_{r\d}+G_{r\d})(\l\vee1)^{-t/2}\\
     \max\bigl\{1,\sqrt{|\log(K_{r\d}+G_{r\d})|},\sqrt{\log(\l+1)}\bigr\};
   \end{multline}
if $K_{r\d}+G_{r\d} \ge 1$ and $\l^{(r+1)/2} \ge K_{r\d}+G_{r\d}$, then there is a 
constant~$\a_{5t}$ such that
\eq\label{Th2-2}
   \|P_X-\n\| \Le \a_{5t}(K_{r\d}+G_{r\d})\l^{-t/2}\max\bigl\{1,\sqrt{\log(\l+1)}\bigr\}.
\en
\end{theorem}

\medskip
\proof
For $K_{r\d}+G_{r\d} < 1$ and $\l\ge1$, we use both parts of~\Ref{tails-of-nu},
with 
\[
    m_1 \Def \lfloor \l -  c_{r\l}\sqrt{\l\log(\l+1)} \rfloor \quad\mbox{and}\quad 
   m_2 \Def \lceil\l + r +  c_{r\l}\sqrt{\l\log(\l+1)} \rceil,
\]
where $\lfloor x \rfloor \le x \le \lceil x \rceil$ denote the integers closest
to~$x$, obtaining
\[
   |\n|\{[0, m_1]\cup[m_2,\infty)\} \Le 2B_r(\l+1)^{-c_{r\l}^2/3} 
      \Le 2B_r(K_{r\d}+G_{r\d})(\l+1)^{-(r+1)/2},
\]
if $c^2_{r\l} := 3(r+1)/2 + |\log(K_{r\d}+G_{r\d})|/\log(\l+1)$.
Hence, from \Ref{TV-1}--\Ref{TV-2}, it follows that
\eqs
   \|P_X-\n\| &\le& \{2c_{r\l}\sqrt{\l\log{\l+1)}} +r+ 2\}\a'_{1t} 
         (K_{r\d}+G_{r\d}) \l^{-(t+1)/2}\\
   &&\mbox{}\qquad + 2\a'_{2t} (K_{r\d}+G_{r\d}) \l^{-t/2}  
         + 4B_r(K_{r\d}+G_{r\d})\l^{-(r+1)/2},
\ens
so that 
\begin{multline*}
   \|P_X-\n\| \Le \b_{3t}(K_{r\d}+G_{r\d})\l^{-t/2}\\
     \max\bigl\{1,\sqrt{\log(1/(K_{r\d}+G_{r\d}))},\sqrt{\log(\l+1)}\bigr\},
   \end{multline*}
with $\b_{3t} := \a'_{1t}\{\sqrt{6(r+1)} + r + 4\} + 2\a'_{2t} + 4B_r$.

For $K_{r\d}+G_{r\d} < 1$ and $\l<1$, we take 
$m_2 := \bigl\lceil \l + r + \sqrt{3|\log(K_{r\d}+G_{r\d})|}\bigr\rceil$
in~\Ref{tails-of-nu}, giving
\[
   |\n|\{[m_2,\infty)\} \Le B_r(K_{r\d}+G_{r\d}),
\]
and
\eqs
  \|P_X-\n\| &\le& (r+2+ \sqrt{3|\log(K_{r\d}+G_{r\d})|})\a'_{1t} (K_{r\d}+G_{r\d})\\
    &&\mbox{}\qquad + 2\a'_{2t} (K_{r\d}+G_{r\d}) + 2B_r(K_{r\d}+G_{r\d}),
\ens
so that 
\[
  \|P_X-\n\| \Le \b'_{3t}(K_{r\d}+G_{r\d})
     \max\bigl\{1,\sqrt{|\log(K_{r\d}+G_{r\d})|},\sqrt{\log(\l+1)}\bigr\},
\]
with $\b'_{3t} := \a'_{1t}\{r+4\} + \a_{2t} + 2B_r$.
Then~\Ref{Th2-1} follows, with $\a_{4t} := \max\{\b_{3t},\b'_{3t}\}$.

For $\l^{t/2} \ge K_{r\d}+G_{r\d} \ge1$, we take $m_1 := \lfloor \l - c_r\sqrt{\l\log(\l+1)} \rfloor$
and $m_2 := \bigl\lceil \l + r + c_r\sqrt{\l\log(\l+1)} \bigr\rceil$, with
$c_r := \sqrt{3t/2}$, giving
\[
   |\n|\{[0, m_1]\cup[m_2,\infty)\} \Le 2B_r(\l+1)^{-t/2} 
\]
Using \Ref{TV-1} and~\Ref{TV-2}, it follows that
\eqs
    \|P_X-\n\| &\le& \{2c_{r}\sqrt{\l\log(\l+1)} +r+ 2\}\a'_{1t} (K_{r\d}+G_{r\d})
         \l^{-(t+1)/2}\\
    &&\mbox{}\qquad+ 2\a'_{2t} (K_{r\d}+G_{r\d}) \l^{-t/2} + 4B_r\l^{-t/2}\\
   &\le& \a_{5t}(K_{r\d}+G_{r\d})\l^{-t/2}\max\bigl\{1,\sqrt{\log(\l+1)}\bigr\},
\ens
with $\a_{5t} = \a'_{1t}(\sqrt{6(r+1)}+r+2) + 2\a'_{2t} + 4B_r$.

\medskip
Note that if $K_{r\d}+G_{r\d} \ge 1$ and $\l^{t/2} < K_{r\d}+G_{r\d}$, one cannot hope to 
get a useful bound from Theorem~\ref{Th1}. If $\l \ge 1$, the error bound for the
individual probabilities is then at least of size~$\l^{-1/2}$, which is the
same size as many of the probabilities themselves. If $\l < 1$, the bound on the individual
probabilities is of size comparable to~$1$.
\ep

\remark  Under the extra conditions that $\ps$ is twice differentiable
and that either \Ref{local-bnd} or~\Ref{local-bnd-tilde} holds with $r\ge2$, Proposition~\ref{Th0-tv} shows that the
factor $\sqrt{\log(\l+1)}$ can in fact be dispensed with.  Note that, to
satisfy the conditions of the proposition, it is necessary to take
$\ch(\th) := \exp\{e^{i\th}-1-i\th\}$, to get $u'(0)=0$.

\medskip
Sometimes it is convenient, for simplicity, to use parameters in the
expansions that are not those emerging naturally from the proofs.  The
following theorem shows that such alterations can easily be allowed for.

\medskip
\begin{theorem}\label{new-pars}
   Suppose that
\[
     \f_\m \Def p_\l A;\qquad \f_{\n\ui} \Def p_\l A';\qquad
     \f_{\n\ut} \Def p_{\l'}A,
\]
with $A(\th) := 1 + \slir a_l \th^l$, $A'(\th) := 1 + \slir a'_l \th^l$
and with $\l > \l'$.  Then, with $\r := 2\pi^{-2}\l$, $\r' := 2\pi^{-2}\l'$
and $a_0 := 1$,
\eqs
     \dloc(\m,\n\ui) &\le& \slir \a_{1l}|a_l - a'_l| (\r\vee1)^{-(l+1)/2};\\
     \dK(\m,\n\ui) &\le&  \slir \a_{2l}|a_l - a'_l| (\r\vee1)^{-l/2};\\
     \dloc(\m,\n\ut) 
        &\le& (\l-\l') \sum_{l=1}^{r+1} \a_{1l} |a_{l-1}|(\r'\vee1)^{-(l+1)/2};\\
     \dK(\m,\n\ut) &\le& (\l-\l') \sum_{l=1}^{r+1} \a_{2l} |a_{l-1}|(\r'\vee1)^{-l/2}.
\ens
\end{theorem}
 
\proof
   For the comparison between $\m$ and~$\n\ui$, we have
\[
      |A(\th) - A'(\th)| \Le \slir |a_l - a'_l| \,|\th|^l, \quad 0 < |\th| \le \pi,
\]
and Proposition~\ref{Th0'} completes the proof.  For that between $\m$ and~$\n\ut$,
note that $p_\l = p_{\l-\l'}p_{\l'}$, and that, for $\l>\l'$ and $0 < |\th| \le \pi$, 
\[
      |p_{\l-\l'}(\th) - 1|\,|A(\th)| 
         \Le (\l-\l')|\th|\,\Blb 1 + \slir |a_l||\th|^l \Brb,
\]
from which and Proposition~\ref{Th0'} the remaining results follow.
\ep

\setcounter{equation}{0}
\section{Poisson approximation}\label{Poisson}
The measures~$\n_r$ considered above are very explicit.
Nevertheless, it is even neater to have approximation in terms
of a Poisson distribution, where possible.  Clearly, if~\Ref{local-bnd}
holds for any $r=r_0$, $\d=\d_0$, then it holds with $r=0$ 
and~$\ps_0(\th)=1$ for all~$\th$, with the exponent $r+\d$ replaced 
by~$\d_0$ if $r_0=0$ and by~$1$ if $r_0 \ge 1$, with~$K_0$
depending on $K_{r_0}$ and on $a_1,\ldots,a_{r_0}$.
Theorem~\ref{Th1} then gives approximation by~$\Po(\l)$ with
accuracy in Kolmogorov distance of order $O(\l^{-t_0/2})$, for
$t_0 = \min\{1,r_0+\d_0\}$.  

However, if $r_0 \ge 1$, one can also write
\[
    \ps(\th)\pl(\th) \Eq \hps(\th)p_{\l'}(\th),
\]
for any $\l' > 0$, where
\[
      \hps(\th) \ :=\ \ps(\th)\exp\{(\l-\l')(e^{i\th}-1)\}.
\]
Taking $\l'-\l=a_1$ now gives a bound
\[
     |\hps(\th) - 1| \Le K'_1 |\th|^{t_1},
\]
of the form~\Ref{local-bnd}, with $t_1 = \min\{r_0+\d_0,2\}$.  
Hence, Theorem~\ref{Th1} implies the following approximation.

\medskip
\begin{corollary}\label{mean}
    If~$X$ has characteristic function $\f_X(\th) = \ps(\th)\pl(\th)$
such that~\Ref{local-bnd} is satisfied with $r \ge 1$, then we have
\eqs
  1.&& \dloc(P_X,\Po(\l')) \Le \a_{1t}  \kkk
     (\r'\vee1)^{-(t+1)/2};\\
  2.&& \dK(P_X,\Po(\l')) \Le \a_{2t}  \kkk 
     (\r'\vee1)^{-t/2},
\ens
where $\l'=\l+a_1$, \hbox{$t := \min\{2,r+\d\}$} and $\r' = 2\p^{-2}\l'$.  
\end{corollary}

\medskip
\nin  The parameter~$\l'$
is chosen to make the Poisson mean~$\l'$ equal to the mean $\l+a_1$ of~$X$.
This choice of the Poisson parameter
improves the rate, in the asymptotic sense that, if $a_1,\ldots,a_r$ 
and~$K_{r\d}$ remain bounded but $\l\to\infty$,  and if $r+\d \ge 2$,
then the approximation error for Kolmogorov distance
is of order~$O(\l^{-1})$, as opposed to the rate of order~$O(\l^{-1/2})$
in general obtained when approximating by~$\Po(\l)$. 

Analogously, fitting the second moment as well (if it is finite)
can lead to further
improvement.  If Poisson approximation is still the aim, the
easiest way to proceed is to consider translating the random
variable~$X$, and approximating $X-m$ by a Poisson instead: 
now one would wish to fix
\[
    \l' \Eq \var X \Eq \ex X - m.
\]
This works well if $\langle \var X - \ex X \rangle = 0$,
where $\langle x \rangle$ denotes the fractional part of~$x$,
but fails otherwise, since $X-m$ only remains integer valued
if~$m$ is itself an integer.  For general~$X$, we therefore
use an average of two adjacent Poisson probabilities
to approximate $P_X\{j\}$. The details are as follows.

Suppose that~\Ref{local-bnd} is satisfied with $r \ge 2$:
$\f_X(\th) = \ps(\th)\pl(\th)$ with
\[
    |\ps(\th) - \ps_r(\th)| \Le K_{r\d} |\th|^{r+\d},
\]
where $\ps_r(\th) = \sjor a_j(i\th)^j$.
For $m\in\integ$ and $0\le p \le 1$,
define the probability measure $ Q_{\l'mp}$ by
\eq\label{Q-def}
  Q_{\l'mp}\{j\} \ :=\ p\Po(\l')\{j-m-1\} + (1-p)\Po(\l')\{j-m\},
\en
having characteristic function $q_{\l'mp}$ given by
\eq\label{q-def}
  q_{\l'mp}(\th) \ :=\ e^{im\th}(1+p(e^{i\th}-1))p_{\l'}(\th).
\en
For $p=0$, $Q$ has the distribution of~$Z'+m$, where $Z'\sim\Po(\l')$;
for $p=1$, $Q$ has the distribution of~$Z'+m+1$; for $0<p<1$,
$Q$ is a mixture of these two distributions. Thus the family of 
distributions
$Q_{\l' mp}$ can be interpreted as a natural generalization of
the usual translated Poisson family, in which the translation
is not restricted to the integers, but may take any real value.
Then we can equivalently write $\f_X(\th) \Eq \bps(\th)q_{\l'mp}(\th)$,
with
\eq\label{bps-def}
   \bps(\th) \Def \ps(\th)\exp\{(\l-\l')(e^{i\th}-1)-im\th\}\{1+p(e^{i\th}-1)\}^{-1}.
\en
In the expansion of~$\bps(\th)$, the coefficients of $\th$ and~$\th^2$
are equal to zero if
\[
     \ex X \Eq \l+a_1 \Eq m+\l'+p \andy 
        \var X \Eq \l + (2a_2 - a_1^2) \Eq \l' + p(1-p).
\]
Then $m\in\integ$, $0\le p < 1$ and~$\l'$ satisfy these two equations
if
\eq\label{mpl-def}
\begin{split}
  m &\Def \lfloor a_1 - (2a_2-a_1^2) \rfloor;\quad
  p^2 \Def \langle a_1 - (2a_2-a_1^2) \rangle;\\
  \l' &\Def \l + (2a_2-a_1^2) - p(1-p),
\end{split}
\en
and it then follows that
\eq\label{kappa-def}
    |\bps(\th) - 1| \Le  \kkk|\th|^t,
\en
for suitable choice of $ \kkk$ depending on $a_1\ldots a_r$ and~$K_{r\d}$,
with $t = \min\{3,r+\d\}$.  Since also, from~\Ref{q-def} and~\Ref{rho-def},
\[
     |q_{\l'mp}(\th)| \Le |p_{\l'}(\th)| \Le e^{-\r'\th^2},
\]
with $\r' = 2\p^{-2}\l'$, the conditions of Proposition~\ref{Th0}
are satisfied with $\ch = q_{\l'mp}$, yielding the following corollary.

\medskip
\begin{corollary}\label{mean-variance}
    If~$X$ has characteristic function $\f_X(\th) = \ps(\th)\pl(\th)$
such that~\Ref{local-bnd} is satisfied with $r \ge 2$, then, for
$\l'$, $m$ and~$p$ defined as in~\Ref{mpl-def} and for $t := \min\{3,r+\d\}$,
we have translated Poisson approximation of the form
\eqs
  1.&& \dloc(P_X, Q_{\l'mp}) \Le \a_{1t}  \kkk
     (\r'\vee1)^{-(t+1)/2};\\
  2.&& \dK(P_X, Q_{\l'mp}) \Le \a_{2t}  \kkk 
     (\r'\vee1)^{-t/2},
\ens
where $\r' = 2\p^{-2}\l'$ and~$ \kkk$ is as in~\Ref{kappa-def}.  
If~\Ref{local-bnd} is replaced by~\Ref{local-bnd-tilde}, then one
takes $a_1 := \ta_1$ and $a_2 := \ta_2 + \half\ta_1$ in~\Ref{mpl-def},
to determine $\l'$, $m$ and~$p$.
\end{corollary}
 
\medskip
\nin In particular, if $a_1,\ldots,a_r$ and~$K_{r\d}$ remain bounded but
$\l\to\infty$, and if $r+\d\ge3$, then $t=3$ and the order of approximation
in Kolmogorov distance is of order $O(\l^{-3/2})$.

\setcounter{equation}{0}
\section{More general expansions}\label{CP}

We now consider cases in which the role of the Poisson family~$\Po(\l)$ is  replaced by
that of another family of probability distributions~$R_\l$, $\l\ge1$, on the integers.
We shall assume that, for $Z_\l\sim R_\l$, $\m(\l) := \ex Z_\l$ 
and~$\s^2(\l) := \var Z_\l$ exist, and are both continuous functions of~$\l$,
with~$\s^2(\l)$ increasing to infinity with~$\l$.  Suppose also that there
exist $c>0$ and~$h(\l)$ such that
\eq\label{h-assn}
  \inf_{\l\ge1}\inf_{0 < |\th| \le \p} \frac1{\th^2 h(\l)}
     \bigl\{-\log|r_\l(\th)|\bigr\} \ \ge\ c ,
\en
where~$r_\l$ is the characteristic function of~$R_\l$.  Clearly, if~\Ref{h-assn}
is satisfied, one could take 
\[
     h(\l) \Def \inf_{0 < |\th| \le \p} \frac1{\th^2}
         \bigl\{-\log|r_\l(\th)|\bigr\}
\]
and $c=1$, or else maybe $h(\l) := \s^2(\l)$ with~$c$ to be determined, but 
it may also be more convenient to choose 
some other, simpler form.  Then, much as in Section~\ref{PC}, we can consider
approximating the distribution of a random variable~$X$ with characteristic
function $\f_X := \ps(\th)r_\l(\th)$ by that of a signed measure~$\n_r
= \n_r(R_\l;\ta_1,\ldots,\ta_r)$ with characteristic function 
\[
    \f_{\n_r}(\th) \Def \tps_r(\th)r_\l(\th) 
        \Def \slor \ta_l\bigl(e^{i\th}-1\bigr)^l r_\l(\th),
\]
(as usual, $\ta_0 = 1$). As in the Poisson case, $\n_r$ is just the linear 
combination $\slor \ta_l D^l R_\l$ of the differences $D^l R_\l$ of the 
probability measure~$R_\l$.  Approximation of the characteristic functions
could be expressed either as 
\eq\label{local-bnd-tilde-Gen}
   |\ps(\th) - \tps_r(\th)| \Le K_{r\d} |\th|^{r+\d},\qquad |\th| \le \pi,
\en	 
for real coefficients~$\ta_l$ and for $r\in\nat_0$, $0<\d\le1$, or as 
\eq\label{local-bnd-Gen}
   |\ps(\th) - \ps_r(\th)| \Le K_{r\d} |\th|^{r+\d},\qquad |\th| \le \pi,
\en	 
where~$\ps_r(\th)$ is as in~\Ref{psr-def}, in which case the corresponding
coefficients~$\ta_l$ can be deduced from~\Ref{b-def}.  These considerations
lead to the following theorem, following directly from  Proposition~\ref{Th0}.

\medskip
\begin{theorem}\label{Th1-G}
Let $X$ be a random variable on~$\integ$ with distribution~$P_X$.  Suppose that 
its characteristic function $\f_X$ is of the form~$\ps R_\l$, where 
$R_\l$ is as above. Suppose also 
that~\Ref{local-bnd-tilde-Gen} is satisfied, for some $r\in\nat_0$ and $\d\ge0$. 
Then, writing 
$t=r+\d$, we have 
\eqs
  1.&& \dloc(P_X,\n_r) \Le \a_{1t} K_{r\d}(\r\vee1)^{-(t+1)/2};\\
  2.&& \dK(P_X,\n_r) \Le \a_{2t}  K_{r\d} (\r\vee1)^{-t/2},
\ens
with $\r := c\,h(\l)$, $\a_{1t}$ and $\a_{2t}$ as in Proposition~\ref{Th0}, and 
$$
\n_r = \n_r(R_\l;\ta_1,\ldots,\ta_r).
$$
If~\Ref{local-bnd-tilde-Gen} is replaced by~\Ref{local-bnd-Gen}, the
corresponding bounds hold with $K_{r\d}$ replaced
by~$K_{r\d}+G_{r\d}$, with $G_{r\d}:=\G_r\p^{1-\d}$ and~$\G_r$ as
in~\Ref{Gr}.
\end{theorem}

\medskip
As in Section~\ref{Poisson}, one may prefer to approximate with a
suitably translated member of the family~$\{R_\l,\,\l\ge1\}$, rather
than with a signed measure~$\n_r$. The corresponding
family of distributions $ Q_{mp}(R_\l)$, for $m\in\integ$ and $0\le p \le 1$,
is given by
\eq\label{Q-def-Gen}
  Q_{mp}(R_\l)\{j\} \ :=\ pR_\l\{j-m-1\} + (1-p)R_\l\{j-m\},
\en
having characteristic function $q_{mp}^{(R_\l)}$ given by
\eq\label{q-def-Gen}
  q_{mp}^{(R_\l)}(\th) \ :=\ e^{im\th}(1+p(e^{i\th}-1))r_{\l}(\th).
\en
Once again, the trick is to find $\l'$, $m$ and~$p$ so that the mean and variance
of~$X$ and of the distribution~$Q_{mp}(R_\l)$ are matched.  

If~\Ref{local-bnd-Gen} is 
satisfied with $r \ge 2$, matching mean and variance implies that we need
\eqa
  \ex X &=& \m(\l)+a_1 \Eq m + \m(\l') + p; \non\\
        \var X &=& \s^2(\l) + (2a_2 - a_1^2) \Eq \s^2(\l') + p(1-p),
      \label{mpl-reqt-Gen}
\ena
where the coefficients $a_1$ and~$a_2$ are as in~\Ref{psr-def}.
These equations have a solution, as long as $\var X \ge \s^2(1) + 1/4$, obtained 
as follows.
For $0 \le p \le 1$, let $\l(p)$ be defined to be the solution of the
equation $\s^2(\l(p)) = \var X - p(1-p)$, noting that $\l(0) = \l(1)$.
Choose 
$m^* := \lfloor \ex X - \m(\l(0)) \rfloor$.  Then the continuous function 
\[
   f(p) \Def \ex X - \m(\l(p)) - m^* - p
\]
satisfies $f(0) \ge 0 > f(1)$, so that there exists a~$p^*$ such that $f(p^*)=0$.
Then the choice $\l' = \l(0)$, $m^*$ and~$p^*$ satisfies~\Ref{mpl-reqt-Gen}, as
desired.

\medskip
\begin{corollary}\label{mean-variance-Gen}
If~$X$ has characteristic function $\f_X(\th) = \ps(\th)R_\l(\th)$
and if~\Ref{local-bnd-Gen} is satisfied with $r \ge 2$, then, for
$\l'$, $m$ and~$p$ solving~\Ref{mpl-reqt-Gen} and for $t := \min\{3,r+\d\}$,
we have translated $R_\l$-approximation of the form
\eqs
  1.&& \dloc(P_X, Q_{mp}) \Le \a_{1t}  \kkk
     (\r'\vee1)^{-(t+1)/2};\\
  2.&& \dK(P_X, Q_{mp}) \Le \a_{2t}  \kkk 
     (\r'\vee1)^{-t/2},
\ens
where $\r' = c\,h(\l')$ and for suitable choice of~$ \kkk$.  
\end{corollary}

\medskip
The most natural application of the above theorem is to mod-compound
Poisson approximation.  For $\l > 0$ and for~$\m$ a probability
distribution on~$\integ$, let $\CP(\l,\m)$ denote the distribution
of the sum $Y := \sum_{j\in\integ\setminus\{0\}} jZ_j$, where $Z_j$, 
$j\ne0$, are independent, and $Z_j \sim \Po(\l\m_j)$.  Then,
if $\m_1>0$, the
characteristic function of~$Y$ is of the form $R_\l := \zl p_{\l_1}$,  
where~$\zl$ is the characteristic function of 
$\sum_{j\in\integ\setminus\{0,1\}} jZ_j$ and $\l_1 = \l\m_1$.
Thus, for the purposes of applying Theorem~\ref{Th1-G} and 
Corollary~\ref{mean-variance-Gen}, $\r$ can be
taken to be $2\p^{-2}\l_1$.  Corollary~\ref{mean-variance-Gen}, for
instance, then gives conditions under which translated compound 
Poisson distribution can be achieved, with approximation at rate
$O(\l^{-3/2})$. 

These considerations apply as long as $\m_1 > 0$, and could
also be invoked if $\m_{-1} > 0$.  If $\m_1=\m_{-1}=0$,
there is then no factor of the form~$\pl$ to guarantee that, for some
$\r>0$, the characteristic function $\f_Y$ of~$Y$ 
(corresponding to the characteristic function~$\ch$ of Proposition~\ref{Th0})
satisfies $|\f_Y(\th)| \le \exp\{-\r\th^2\}$ for all
$|\th| \le \p$.  Some additional aperiodicity condition needs to be
satisfied, if the family $\{\CP(\l,\m),\,\l\ge1\}$ is to
satisfy~\Ref{h-assn}.  Indeed, if $Y = 2Z$ where $Z\sim\Po(\l)$, 
and if $W\sim\Be(1/2)$ is independent of~$Y$, it is not
true that the distribution of $Y + W$ is close to that of~$Y$ in total
variation, even though $|\f_{Y+W}(\th) - \f_Y(\th)| \le K_0|\th|$.

\setcounter{equation}{0}
\section{Applications}\label{Apps}

\subsection{A single convolution}\label{convolution}
The most obvious application of the above results arises when $\f_X = \ps\pl$
and~$\ps$ is itself the characteristic function of a probability distribution
on the integers.  In this case, $X$ is the sum of two independent random
variables, one of them with the~$\Po(\l)$ distribution, and the situation
is probabilistically very simple.  For example, we could take~$\ps$ to be
the characteristic function of a random variable~$Y_s$ with
\[
    \pr[Y_s=j] \Eq s!\,s\sji \frac1{j(j+1)\ldots(j+s)}
\]
for some integer $s\ge1$.  Calculation shows that~$Y_s$ has characteristic
function
\[
    \ps(\th) \Eq 1 + s\sum_{l=1}^{s-1}\frac{(1-e^{-i\th})^l}{s-l}
        - s(1-e^{-i\th})^s \log(1-e^{i\th}),
\]
and that~\Ref{local-bnd-tilde} holds with $r=s-1$ and any $\d < 1$ if
\[
     \ta_0 \Eq 1;\qquad \ta_j \Eq \sum_{l=1}^j (-1)^{j-l}\,\frac s{s-l}\,
       {j-1\choose l-1}, \quad 1\le j\le s-1.
\]
Hence, if $X = Z + Y_s$, where $Z\sim\Po(\l)$ is independent of~$Y_s$, then
the theorems in Sections \ref{PC} and~\ref{Poisson} can be applied,
provided that~$s$ is large enough; in particular, a translated Poisson
approximation can be applied with accuracy of order~$O(\l^{-3/2+\e})$
for any $\e > 0$ if $s=3$ (in which case~$X$ has finite second moment),
and of order~$O(\l^{-3/2})$ if $s\ge4$.  Similar considerations
apply to the approximation of $X = Z - Y_s$.

\subsection{Sums of independent random variables}

Let $X_1,\ldots,X_n$ be independent integer valued random variables,
and let~$S_n$ denote their sum.  In contexts in which a central limit
approximation to the distribution of~$S_n$ would be appropriate, the
classical Edgeworth expansion (see, e.g., Petrov 1975, Chapter~5) is
unwieldy, because~$S_n$ is confined to the integers.  As an
alternative, Barbour and \ceka~(2002) give a Poisson--Charlier
expansion, for~$S_n$ `centered' so that its mean and variance are
almost equal, with an error bound expressed in the total variation
norm.  Here, we show that such an expansion can be justified by the
techniques of this paper.

   Assume that each of the $X_j$ has finite $(r+1+\d)$'th moment, with
$r\ge1$, and define
\eq\label{B-def}
     A\ur(w) \Def 1 + \slt \ta_l\ur w^l 
       \Eq \exp\Blb \sum_{l=2}^{r+1} \frac{\k_l w^l}{l!}\Brb,
\en
where $\k_l := \k_l(S_n)$ and $\k_l(X)$ denotes the $l$'th
factorial cumulant of the random variable~$X$.  Then the approximation
that we establish is to the Poisson--Charlier signed measure~$\n_r$
with
\eq\label{nu-def-indept}
    \n_r\{j\} \Def \Po(\l)\{j\}\Blb 1 + \sum_{l=2}^{L_r} (-1)^l
       \ta_l\ur C_l(j;\l)\Brb,
\en
where $L_r := \max\{1,3(r-1)\}$, and where $\l := \ex S_n$;
$\n_r$ has characteristic function
\eq\label{nu-char-indept}
   \f_{\n_r} \Def p_\l(\th)\,\tA\ur(\th),
\en
where
\eq\label{tilde-A-def}
  \tA\ur(\th) \Def  1 + \sum_{l=2}^{L_r} \ta_l\ur (e^{i\th}-1)^l .
\en
We need two further quantities involving the~$X_j$:
\eq\label{centering}
   K\un \Def \Bigl|\sum_{j=1}^n \k_2(X_j) \Bigr| 
\en
and  
\eq\label{dtv}
    p_j \Def 1 - \dtv(\law(X_j),\law(X_{j+1})).
\en

\medskip
\begin{theorem}\label{indept}
     Suppose that there are constants $\mmm_l$, $1\le l\le r+1$, such
that, for each~$j$, 
\[
     |\k_l(X_j)| \Le \mmm_l,\quad 2\le l\le r+1;\qquad \ex|X_j|^{r+1+\d}
       \Le \mmm_1^{r+1+\d}.
\]
Suppose also that $p_j \ge p_0 > 0$ for all~$j$, and that $\l\ge n\l_0$.
Then
\[
   \dK(\law(S_n),\n_r) \Le G(\mmm_1,\ldots,\mmm_{r+1},K\un,p_0^{-1},\l_0^{-1})
     n^{-(r-1+\d)/2},
\]
for a function~$G$ that is bounded on compact sets.
\end{theorem}

\medskip \remark For asymptotics in~$n$, with triangular arrays of
variables, the error is of order $O(n^{-(r-1+\d)/2})$ when $\l_0$
and~$p_0$ are bounded away from zero, and $\mmm_1,\ldots,\mmm_{r+1}$
and~$K\un$ remain bounded.  The requirements on $\l_0$ and~$p_0$ can
often be achieved by grouping the random variables appropriately,
though attention then has to be paid to the consequent changes in
the~$\mmm_l$. The final condition can always be satisfied with $K\un
\le 1$, by replacing the~$X_j$ by translates, where necessary.  For
more discussion, we refer to Barbour and \ceka~(2002).  The above
conditions are designed to cover sums of independent random variables,
each of which has non-trivial variance, has uniformly bounded
$(r+1+\d)$'th moment, and whose distribution overlaps with its unit
translate.
 
\medskip
\proof
We check the conditions of Proposition~\ref{Th0'}.  First, in view
of~\Ref{dtv}, we can write
\[
   \ex\bigl(e^{i\th X_j}\bigr) \Eq \half p_j (e^{i\th}+1)\f_{1j}(\th)
      + (1-p_j)\f_{2j}(\th),
\]
where both $\f_{1j}$ and~$\f_{2j}$ are characteristic functions.
Hence we have
\[
   \bigl|\ex\bigl(e^{i\th X_j}\bigr)\bigr|
    \Le 1-p_j + p_j\cos(\th/2) \Le 1 - p_j \th^2/4\pi, \quad 0 \le |\th| \le \pi.
\]
Hence $\f_\m(\th) := \ex\bigl(e^{i\th S_n}\bigr)$ satisfies
\eq\label{mu-bnd-indept}
    |\f_\m(\th)| \Le  \exp\{-np_0\th^2/4\pi\}, \quad 0 \le |\th| \le \pi.
\en
On the other hand, from the additivity of the factorial cumulants, we have
\[
      |\k_l(S_n)| \Le n\mmm_l, \quad 3\le l\le r+1,
\]
with $|\k_2(S_n)| \le K\un$ from~\Ref{centering}. From~\Ref{B-def},
we thus deduce the bound  $|\ta_l\ur| \le c_l n^{\lfloor l/3\rfloor}$,
for $c_l = c_l(K\un,\mmm_3,\ldots,\mmm_{r+1})$, $l\ge1$. Hence
\eq\label{nu-bnd-indept}
   |\f_{\n_r}(\th)| 
     \Le \exp\{-2n\l_0\th^2/\pi^2\}c'n^{\lfloor L_r/3 \rfloor} 
        \Le \exp\{-n\l_0\th^2/\pi^2\}c'',
\en
for $c''=c''(K\un,\mmm_3,\ldots,\mmm_{r+1})$,
and we can take $\h := C e^{-n\r'\th_0^2}$ in
Proposition~\ref{Th0'}, for
$$
\r' = \min\{\l_0/\pi^2,p_0/4\pi\}
$$ 
and a suitable $C = C(K\un,\mmm_3,\ldots,\mmm_{r+1})$.  The choice
of~$\th_0$ we postpone for now.

\medskip
For $|\th| \le \th_0$, we take $\chi(\th) := p_\l(\th)$, and check
the approximation of 
$$
  \f_\m(\th)\exp\{-\l(e^{i\th}-1)\} \Eq \ex\Blb (1 + w)^{S_n} \Brb e^{-w\ex S_n}
$$
by $\tA\ur(\th)$ as a polynomial in $w := e^{i\th}-1$. We begin with the inequality
\eqs
   \Bigl|(1 + w)^{s}  - \sum_{l=0}^{r+1} \frac{w^l}{l!}\, s_{(l)}\Bigr|
   &\le&  \frac{|s_{(r+2)}|}{(r+2)!}\,|w|^{r+2} \wedge 2\frac{|s_{(r+1)}|}{(r+1)!}\,|w|^{r+1}\\
   &\le& \frac{|s_{(r+1)}|}{(r+2)!}\,|w|^{r+1+\d}\{|s|+r+1\}^\d \{2(r+2)\}^{1-\d},
\ens
derived using Taylor's expansion, true for any $s\in\integ$ and $0 < \d \le 1$,
where $s_{(l)} := s(s-1)\ldots(s-l+1)$.
Hence, for each~$j$, we have
\eq\label{factorial-bnd}
   \Bigl|\ex\Blb (1 + w)^{X_j} \Brb 
     - \sum_{l=0}^{r+1} \frac{\ex\{(X_j)_{(l)}\}}{l!}\, w^l \Bigr|
          \Le c_{r,\d} |\th|^{r+1+\d}(\mmm_1 + \mmm_1^{r+1+\d}),
\en
for a universal constant $c_{r,\d}$.  Then,
writing
\[
   Q^{(s)}_{r+1}(w;X) \Def \exp\Blb \sum_{l=s}^{r+1}\k_l(X)w^l/l! \Brb,
\]
and using the differentiation formula in Petrov~(1975, p.~170), we have
\begin{align}
   \Bigl|Q\ui_{r+1}(w;X_j) - \sum_{l=0}^{r+1} w^l \ex {X_j \choose l}\Bigr|
   &\le \frac{|\th|^{r+2}}{(r+2)!}\,\sup_{|\th'| \le \th_0}
      \Bigl|\frac{d^{r+2}}{dz^{r+2}} Q\ui_{r+1}(z;X_j)
      \Bigr|_{z=e^{i\th'}-1}
\nonumber\\
   &\le  |\th|^{r+2} c(\mmm_1,\ldots,\mmm_{r+1}), \label{Petrov-1}
 \end{align}
for a suitable function~$c$ and for all~$|\th| \le \pi$.  Combining these 
estimates, we deduce that, for $w = e^{i\th}-1$ and for all $|\th| \le \pi$,
\eq\label{single-diff}
   \Bigl|\ex\Blb (1 + w)^{X_j} \Brb e^{-\ex X_j w} 
      - Q\ut_{r+1}(w;X_j)\Bigr| \Le k_1|\th|^{r+1+\d} ,
\en
where $k_1 = k_1(\mmm_1,\ldots,\mmm_{r+1})$.

Now a standard inequality shows that, for $u_j := \prod_{l=1}^j x_l
\prod_{l=j+1}^n y_l$, for complex $x_l,y_l$ with $y_l \neq 0$ and 
$|x_l/y_l - 1| \le \e_l$, then
\eq\label{u-bnd}
      |u_n-u_0| \Le  |u_0| \Blb \prod_{s=1}^{n-1}(1+\e_s) \Brb \sum_{l=1}^n \e_l .
\en
Taking $x_j :=  \ex\Blb (1 + w)^{X_j} \Brb e^{-\ex X_j w}$ and
$y_j := Q\ut_{r+1}(w;X_j)$, \Ref{single-diff} shows that we can take 
$\e_l := \e := k_1|\th|^{r+1+\d}e^M$ for each~$l$, with 
$$
M := \exp\{\sum_{l=2}^{r+1}\mmm_l/l!\},
$$
provided that $|\th| \le \th_0 \le 1$. Choosing $\th_0 := n^{-1/3}$ then
ensures that $(1+\e)^n$ is suitably bounded, and~\Ref{u-bnd} yields
\eq\label{n-diff}
    \Bigl|\ex\Blb (1 + w)^{S_n} \Brb e^{-w\ex S_n} -  Q\ut_{r+1}(w;S_n)\Bigr|
     \Le k_2 n|\th|^{r+1+\d},
\en
for $k_2 = k_2(K\un,\mmm_1,\ldots,\mmm_{r+1})$, since
\[
   |Q\ut_{r+1}(w;S_n)| \Le \exp\{|\k_2(S_n)|\th_0^2/2\}
      \exp\Blb \sum_{l=3}^{r+1}n\mmm_l |\th_0|^l/l! \Brb
\]
is bounded for $\th_0 = n^{-1/3}$, in view of~\Ref{centering}.

The remaining step is to note that, for $w=e^{i\th}-1$,
\eq\label{Petrov-2}
   \Bigl|Q\ut_{r+1}(w;S_n) - \tA\ur(\th)\Bigr|
    \Le \frac{|\th|^{L_r+1}}{(L_r+1)!}\,\sup_{|\th'| \le \th_0}
      \Bigl|\frac{d^{L_r+1}}{dz^{L_r+1}} Q\ut_{r+1}(z;S_n) \Bigr|_{z=e^{i\th'}-1} ,
\en
where the right hand side is at most $k_3 n^{r-1}|\th|^{L_r+1}(1 + n|\th|^2)$
in $|\th| \le n^{-1/3}$,
with $k_3 = k_3(K\un,\mmm_1,\ldots,\mmm_{r+1})$. Here, we use the facts that 
$|\k_2(S_n)|$ is bounded by~$K\un$, and that each~$\k_l(S_n)$ for $l\ge3$, for which we have
only the weak bound~$n\mmm_l$, occurs associated with the power $w^l$ in the exponent 
of $Q\ut_{r+1}(w;S_n)$.  Combining this with~\Ref{n-diff}, we have
established that for $|\th|\leq n^{-1/3}$, we have
\eq\label{part-1-bnd}
  |\f_\m(\th)\exp\{-\l(e^{i\th}-1)\} - \tA\ur(\th)|
   \Le k_4 n|\th|^{r+1+\d}(1 + (n|\th|^2)^{r-1}), 
\en
where  $k_4 = k_4(K\un,\mmm_1,\ldots,\mmm_{r+1})$.  This gives 
\begin{gather*}
\kkk_1 = nk_4,\quad
t_1=r+1+\d,\quad 
\kkk_2 = n^{r}k_4,\quad
t_2 = 3r-1 + \d\\
\kkk = ,\quad \r = 2\l/\pi^2,\quad \e=0,\quad\text{ and } \th_0 =
n^{-1/3}
\end{gather*}
in Proposition~\ref{Th0'}, together with $\h = Ce^{-n^{1/3}\r'}$ from
the earlier bounds.  Applying Corollary~\ref{Cor0}, and using the tail
properties of the Poisson--Charlier measures~\Ref{tails-of-nu}, the
theorem follows.  \ep

\nin A total variation bound of precisely the same order can also be deduced, 
by combining the arguments used for Propositions \ref{Th0'} and~\ref{Th0-tv}.
Note that~$\f_\m$ is twice differentiable, because the~$X_j$ all have finite
second moments, and that, as in Section~\ref{PC}, we need to take 
$\ch(\th) := \exp\{e^{i\th}-1-i\th\}$.

\subsection{Analytic combinatorial schemes}\label{acs}
An extremely interesting ran\-ge of applications is to be found in the
paper of Hwang~(1999).  His conditions are motivated
by examples from combinatorics, in which generating functions are
natural tools.  He works in an asymptotic setting, assuming that~$X_n$
is a random variable whose \pgf~$R_n$ is of the form 
\[
    R_n(z) \Eq z^h(g(z) + \e_n(z))e^{\l(z-1)},
\]
where $h$ is a non-negative integer, 
and both $g$ and~$\e_n$ are analytic in a closed disc of radius $\h > 1$.
As $\nti$, he assumes that $\l\to\infty$ and that
$\sup_{z:|z|\le\h}|\e_n(z)| \le K\l^{-1}$, uniformly in~$n$.  He then
proves a number of results describing the accuracy of the approximation
of~$P_{X_n-h}$ by $\Po(\l + g'(1))$.

Under his conditions, it is immediate that we can write
\eq\label{ge-series}
    g(z) \Eq \sjo g_j(z-1)^j \andy \e_n(z) \Eq \sjo \e_{nj}(z-1)^j
\en
for $|z| < \h-1$, with 
\eq\label{ge-bnds}
    |g_j| \Le k_g (\h-1)^{-j} \andy |\e_{nj}| \Le \l^{-1}k_\e (\h-1)^{-j}
\en
for all $j\ge0$.  Hence $X := X_n-h$ has characteristic function of
the form $\ps\pl$, where 
\[
    \ps\un(\th) \Eq g(e^{i\th}) + \e_n(e^{i\th}),
\]
and hence, for any $r\in\nat_0$,
\eq\label{appx-1}
     |\ps\un(\th) - \tps\un_r(\th)| \Le K_{r1} |\th|^{r+1},\qquad |\th| \le (\h-1)/2,
\en
with~$\tps$ defined as in~\Ref{psr-def-tilde}, taking $\ta\un_j = g_j+\e_{nj}$;
note that the constant~$K_{r1}$ can indeed be taken to be uniform for all~$n$.
Since also $g$ and $\e_n$ are both uniformly bounded on the unit circle,
and since~$\tps_n$ is bounded (uniformly in~$n$) for~$|\th| \le \p$, it is
clear that~\Ref{appx-1} can be extended to all $|\th| \le\p$, albeit with
a different uniform constant~$K'_{r1}$, so that~\Ref{local-bnd-tilde} holds
with $\d=1$ for any $r\in\nat_0$.  Thus Theorems \ref{Th1} and~\ref{Th2}
can be applied with any choice of~$r$, giving progressively more
accurate approximations to $P_{X_n-h}$, as far as the $\l$-order is
concerned, in terms of progressively more 
complicated perturbations of the Poisson distribution.  These
theorems are thus applicable to all the examples that Hwang
considers, including the numbers of components (counted in
various ways) in a wide class of logarithmic assemblies, 
multisets and selections.

For instance,
Corollary~\ref{mean-variance} gives an approximation to~$P_{X_n-h}$
by the mixture~$Q_{\l'mp}$ with
\[
    m \Def \lfloor m_n - v_n \rfloor;\qquad
  p^2 \Def \langle m_n-v_n \rangle;\qquad
  \l' \Def \l + v_n - p(1-p),
\]
where $m_n := g_n'(1)$, $v_n := g_n''(1) + g_n'(1) - \{g_n'(1)\}^2$
and $g_n := g + \e_n$.  Hwang's approximation by~$\Po(\l + g'(1))$ has 
asymptotically the same mean as ours (and as that of~$X_n-h$), but
a variance asymptotically differing by $\k := g''(1) - \{g'(1)\}^2$
(together with an element arising from $p(1-p)$ which is not
in general asymptotically negligible).  As a consequence, Hwang's
approximation has an error of larger asymptotic order, 
in which the quantity~$\k$ appears; for instance, for Kolmogorov distance,
his Theorem~1 gives an error of order $O(\l^{-1})$, whereas that from
Corollary~\ref{mean-variance} is of order~$O(\l^{-3/2})$.

Although our Poisson expansion theorems are automatically 
applicable under Hwang's conditions, they also apply to examples
that do not satisfy his conditions: that of Section~\ref{convolution}
is one such. Conversely, Hwang's Theorem~2, which establishes
Poisson approximation in the lower tail with good {\it relative\/}
accuracy, cannot be proved using only our conditions; the conclusion 
would not be true, for instance, for the random variable $X-Y_s$ of
Section~\ref{convolution}.

Note also that Hwang examines problems from combinatorial settings
in which approximation is not by Poisson distributions:  he has
examples concerning the Bessel family,
\[
     B(\l)\{j\} \Def L(\l)^{-1}\frac{\l^j}{j!(j-1)!},\quad j\in\nat,
\]
for the appropriate choice of~$L(\l)$.  Here, we could apply
Corollary~\ref{mean-variance-Gen} to obtain slightly sharper
approximations than his within the translated Bessel family,
or Theorem~\ref{Th1-G} to obtain asymptotically more accurate 
expansions.

\subsection{Prime divisors}

The numbers of prime divisors of a positive integer~$n$, counted
either with ($\Om(n)$) or without ($\om(n)$) multiplicity, can also be
treated by these methods, since excellent information is available
about their generating functions.  For our purposes, we use only the
shortest expansion, taken from Tenenbaum~(1995, Theorems II.6.1 and~6.2).
One finds that for~$N_n$ uniformly distributed on $\{1,2,\ldots,n\}$
we have 
\eqs 
\ex\{e^{i\th\om(N_n)}\}
&=& p_{\loglog n}(\th)\Blb \Phi_1(e^{i\th}-1)+ \h_1(\th)\Brb;\\
\ex\{e^{i\th\Om(N_n)}\} 
&=& p_{\loglog n}(\th)\Blb \Phi_2(e^{i\th}-1)
+ \h_2(\th)\Brb, 
\ens 
where $|\h_s(\th)| \le C_s/\log n$, $s=1,2$, for
some constants $C_1$ and~$C_2$, and 
\eqs 
\Phi_1(w)= &:=&
\frac{1}{\Gamma(1+w)}
\prod_q \Bl 1 + \frac w{q}\Br \,\Bl 1 - \frac1q \Br^w;\\
\Phi_2(w) &:=& \frac{1}{\Gamma(1+w)} \prod_q \Bl 1 - \frac w{q-1}
\Br^{-1} \, \Bl 1 - \frac1q \Br^w, 
\ens 
$q$ running here over prime numbers.  These expansions were
established and used by R\'enyi and Tur\'an~(1958) in their proof of
the Erd\H os--Kac Theorem, but they are also sketched by Selberg
(1954). We refer to Kowalski and Nikeghbali (2009) for the structural
interpretation of the two factors in these functions (with
$1/\Gamma(1+w)$ being related to the number of cycles of large random
permutations).

Let $\ta_{ls}$, $s=1,2$, denote the Taylor coefficients of the
functions $\Phi_s(w)$ as power series in~$w$ (around $w=0$, which
corresponds to $\th=0$).  By analyticity, it follows that for any $r$,
we have
\[
   \Blm \Phi_s(w) - 1 - \slir \ta_{ls} w^l \Brm \Le C_{rs}|w|^{r+1},
\]
for suitable constants $C_{rs}$ and for $|w|\leq 2$.  Defining the
measures $\n_r\us$ by
\[
     \n_r\us\{j\} \Def \Po(\loglog n)\{j\}\Bl 1 + \slir (-1)^l \ta_{ls} C_l(j;\loglog n) \Br,
\]
this leads to the following conclusion, which is deduced immediately
from Theorem~\ref{Th1}, and refines the Erd\H os--Kac theorem. 

\begin{theorem}\label{divisors}
For the measures~$\n_r\us$ defined above, we have
\eqs
   \dloc(P_{\om(N_n)},\n_r\ui) 
      &\le& \a'_{1,r+1} C_{r1}(\loglog n)^{-1-r/2} + \ta_1 C_1/\log n;\\
   \dK(P_{\om(N_n)},\n_r\ui) &\le& \a'_{2,r+1} C_{r1}(\loglog n)^{-(r+1)/2} 
        + \tC_1\loglog n/\log n;\\
   \dloc(P_{\Om(N_n)},\n_r\ut) 
      &\le& \a'_{1,r+1} C_{r2}(\loglog n)^{-1-r/2} + \ta_1 C_2/\log n;\\
   \dK(P_{\Om(N_n)},\n_r\ut) &\le& \a'_{2,r+1} C_{r2}(\loglog n)^{-(r+1)/2} 
        + \tC_2\loglog n/\log n,
\ens
for suitable constants $\tC_1$ and~$\tC_2$.
\end{theorem}

\remark Note that it follows from Theorem~\ref{Th2} that the total
variation distance is in each case also of order $O\bigl\{(\loglog
n)^{-(r+1)/2}\bigr\}$.  This can be deduced by applying the theorem to
the expansion with one more term, and then observing that the extra
term has total variation norm of order $O\bigl\{(\loglog
n)^{-(r+1)/2}\bigr\}$, in view of the observation
following~\Ref{BC6.1}. Alternatively, one could use
Proposition~\ref{Th0-tv}.  As far as we know, total variation
approximation was first considered in this context by Harper~(2009),
who proved a bound with error of size $1/(\log\log n)$ (for a
truncated version of $\om(n)$, counting only prime divisors of size up
to $n^{1/(3(\log\log n)^2)}$), and deduced explicit bounds
in Kolmogorov distance.

\medskip
To indicate what this means in concrete terms for number theory
readers, consider the case of $\om(n)$ for $r=1$. Taylor expansion
gives
$$
  \Phi_1(w) \Eq 1 + B_1w + O(w^2)
$$
as $w\rightarrow 0$, where~$B_1\approx 0.26149721$ is the Mertens
constant, i.e., the real number such that
$$
\sum_{{q\leq x}\atop \text{ $q$ prime}}{\frac{1}{q}}=\log\log x+B_1+o(1),
$$
as $x\rightarrow +\infty$.
\par
In view of the remark above, an application of Theorem~\ref{divisors}
gives 
\begin{align*}
\Bigl|\frac{1}{n}|\{k\leq n\,\mid\, \omega(n)\in A\}| - 
   \n_1\ui\{A\}\Bigr|
     &\Le \half\|P_{\om(N_n)} - \n_1\ui\| \\
&     \Eq O\Bl \frac{1}{\log\log n}\Br,
\end{align*}
for any set $A$ of positive integers, where
$$
   \n_1\ui\{j\} \Eq \Po(\log\log n)\{j\} \Bl 1 
+ B_1\Blb 1 - \frac j{\log\log n}\Brb\Br.
$$
Higher expansions could be computed in much the same way.

Alternatively, a more accurate approximation is available from
Corollary~\ref{mean-variance}, while staying within the realm of
(translated) Poisson distributions.

For this, we compute the expansion of $\Phi_1$ to order $2$, obtaining
(after some calculations) that
$$
  \Phi_1(w) \Eq 1+B_1w+ \ta_2w^2 + O(w^3),\quad\text{ as } w\rightarrow
  0,
$$
where
$$
\ta_2 \Def \frac{B_1^2}{2} - \frac{\pi^2}{12} - \frac{1}{2}
\sum_{q\text{ prime}}{\frac{1}{q^2}}
$$
(use $1/\Gamma(1+w)=1+\gamma w+(\gamma^2-\pi^2/12)w^2+O(w^3)$, as well
as the Mertens identity
$$
\gamma+\sum_{\text{$q$ prime}}{\Bigl(\frac{1}{q}+
\log\Bigl(1-\frac{1}{q}\Bigr)\Bigr)}=
B_1,
$$
and expand every term in the Euler product). This corresponds
to~\Ref{local-bnd-tilde}, since $w=e^{i\th}-1$, and therefore we
have~\Ref{local-bnd} with
$$
a_1=B_1,\quad\quad a_2=\ta_2+\frac{1}{2}B_1=
\frac{B_1+B_1^2}{2} - \frac{\pi^2}{12} - \frac{1}{2}
\sum_{q\text{ prime}}{\frac{1}{q^2}}.
$$

We can then apply Corollary~\ref{mean-variance} to get the translated
Poisson approximation~$Q_{\l'mp}$, with parameters calculated
using~\Ref{mpl-def}.  With
$$
x := B_1-(2a_2-B_1^2)=\frac{\pi^2}{6}+\sum_{q\text{ prime}}
{\frac{1}{q^2}} \approx 2.0971815,
$$
this gives
\eqs
   p &=& \sqrt{\langle x \rangle}\ \approx\ 0.31173945; \qquad m \Eq 2;\\
     \l' &=& \log\log n + B_1 - x - p(1-p)\ \approx\ \log\log n 
- 2.0502422
\ens
Thus for any positive integer $n$ and any set $A$ of positive
integers, we have
\begin{multline*}
\Bigl|\frac{1}{n}|\{k\leq n\,\mid\, \omega(n)\in A\}| - 
   \{p\Po(\l')\{A-3\} + (1-p)\Po(\l')\{A-2\}\}\Bigr| \\
     \Eq O\Bl \frac{1}{(\log\log n)^{3/2}}\Br,
   \end{multline*}
   where, again, we can use the total variation norm in view of the
   previous remark.  Similar results hold for $\Omega(n)$, where one
   obtains the following approximate values
 \eqs
    p &\approx&  0.5195; \qquad m \Eq 0;\\
      \l' &\approx&  \log\log n + 0.5152.
 \ens

\end{document}